\documentclass[11pt]{article}
\usepackage{graphicx}
\usepackage{epsfig}
\usepackage{amssymb, amsmath}
\usepackage{amsthm}
\usepackage{setspace}
\usepackage[top = 1in, bottom = 1in, left = 1.2in, right =
1.2in]{geometry}
\newtheorem{theorem}{Theorem}[section]

\newtheorem{lemma}[theorem]{Lemma}

\newtheorem{remark}[theorem]{Remark}

\begin{document}

\title{Derivation of Ohm's law from the kinetic equations}%

\author{Juhi Jang\thanks{Department of Mathematics, University of California Riverside. Email: juhijang@math.ucr.edu Supported in part by NSF Grant DMS-0908007} and Nader Masmoudi\thanks{Courant Institute, New York University. Email: masmoudi@cims.nyu.edu Supported in part by NSF Grant DMS-0703145}}


\singlespacing \maketitle \numberwithin{equation}{section}

\begin{abstract}
The goal  of this article is to 
give a formal  derivation of   Ohm's law of Magnetohydrodynamics (MHD)
 starting  from the Vlasov-Maxwell-Boltzmann system. 
The derivation is based on various physical scalings and the moment methods when the 
Knudsen number goes to zero. We also give a derivation of the so-called Hall 
effect as well as other limit models such as the Navier-Stokes-Maxwell system. 
 Our results include both the compressible and incompressible MHD models. 
\end{abstract}

{\small 2000\textit{\ AMS Subject Classification.} 35Q20, 35Q61, 35Q83, 76W05.{ }{\it  Keywords.} Ohm's law, Vlasov-Maxwell-Boltzmann system,   Magnetohydrodynamics, Hall effect, hydrodynamic limit}

\maketitle

\section{Introduction}

There are different models to describe the state of a  plasma depending on 
several  parameters such as the Debye length, the plasma frequency, the collision 
frequencies between the different species. One of the goals
 of this project is give a 
relation between the several models used, namely how 
to derive rigorously  the simpler models  from the 
more complete ones  and find the regimes where these approximations  are  valid.  
Formal derivation of these models can be found in Plasma Physics text books 
(see for instance  
 Bellan \cite{Bellan06},  Boyd and Sanderson \cite{BS03}, Dendy \cite{Dendy90} and the 
paper \cite{BDDCGT} etc.)
 
Indeed, since the plasma consists of a very large number of interacting particles, 
it is appropriate to adopt a statistical approach. In the kinetic description, 
it is only necessary to evolve the distribution function $f_\alpha (t,x,v)$ for 
each species in the system. Vlasov equation is used in this case with the 
Lorentz force term and collision terms. It is 
 coupled with the  Maxwell equations 
for the electromagnetic  fields.

If collisions are important, then each species is in a local equilibrium 
and the plasma is  treated as a fluid. More precisely it  is treated as a mixture 
of two or more interacting fluids. This is the two-fluid model or the 
so-called  Euler-Maxwell system.  
 
Another level of approximation  consists in treating the plasma as a single fluid 
by using the fact that the mass of the electrons is much smaller than
the mass of the ions or by using that collisions yield that 
all particles will evolve  (at leading order) with the same macroscopic velocity. 
  This  yields  the Hall-MHD (Hall magneto-hydrodynamic  model). 
 Then, one can derive  MHD models if the Hall effect is negligible.

The question on how the kinetic theory and the fluid dynamics are related is very interesting both physically and mathematically and it  
goes back to the founding work of Maxwell \cite{Max} and
Boltzmann \cite{Bol}. Moreover, the purpose of the
Hilbert's sixth problem \cite{Hilbert} is to seek a unified theory
of the gas dynamics including various levels of descriptions from a
mathematical point of view. There has been a lot of important progress on the hydrodynamic limits from the Boltzmann equation over the years; so far there are essentially three different mathematical 
approaches. The first is based on spectral analysis
of the semi-group generated by the linearized Boltzmann equation; see \cite{BU,KMN,Nishida}. The second is based on  Hilbert or
Chapman-Enskog expansions \cite{C,mel,Grad0}; see more recent work in \cite{g3,GJJ09,GJJ10,jj}. The third approach, initiated by Bardos-Golse-Levermore \cite{BGL1,BGL2}, is working in the framework of DiPerna-Lions'
renormalized solutions 
\cite{DL}, to justify global weak solutions of incompressible flows
(Navier-Stokes, Stokes, and Euler), and (compressible) acoustic
system; see \cite{BGL2,BGL3,GL,GS,LM,LM1,LM2,MS,Sa} (see also \cite{Arsenio12} for the 
non-cut-off case).

 However, there are only limited results on hydrodynamic limits for charged particles, which satisfy the Vlasov-Maxwell-Boltzmann (VMB) system, even at the formal level due to the complexity of the system and its underlying  muti-scale feature. One rigorous result is given in \cite{j3} where  a diffusive expansion to VMB system was studied  with one particular scaling in the incompressible regime in the framework of classical solutions \cite{g2}, and as a by-product, new fluid equations interacting with the electric field, where the magnetic effect appears only at a higher order, were derived. See also \cite{belm,GJ10,MT} for hydrodynamic limits from the Boltzmann equations in the presence of an interacting field, but without a magnetic field.
 
  In this article, we are interested in the derivation of MHD type equations, which describe the motion of electrically conducting media in the presence of a magnetic field, from the VMB system by introducing various scalings. As a result, we can also identify the corresponding  Ohm's law of  the MHD equations. 
The goal here is to provide  formal expansions that show the 
relevant scalings. Our analysis is the first step towards a more rigorous analysis. 
 It can be seen as the extension of the paper 
of Bardos-Golse-Levermore \cite{BGL1} to the VMB case.   Some rigorous derivation will 
be given in a  forthcoming paper.

Let us also  mention that various MHD models were obtained from the (macroscopic) two-fluid Euler-Maxwell equations  in \cite{BDDCGT} by taking different asymptotic limits. 
Here, our goal is to start from the kinetic level.

\subsection{Two species Vlasov-Maxwell-Boltzmann system}

The dynamics of charged dilute particles (e.g., electrons and ions)
is described by the Vlasov-Maxwell-Boltzmann system:
\begin{eqnarray}
&&\partial _tF_{+}+v\cdot \nabla _xF_{+}+\frac{e_{+}}{m_{+}}(E+
v \times B)\cdot \nabla _vF_{+}=Q(F_{+},F_{+})+Q(F_{+},F_{-}),
\nonumber
\label{boltzmann} \\
&&\partial _tF_{-}+v\cdot \nabla _xF_{-}-\frac{e_{-}}{m_{-}}(E+
v \times B)\cdot \nabla _vF_{-}=Q(F_{-},F_{+})+Q(F_{-},F_{-}),
\label{VMB} \\
&&F_{\pm }(0,x,v)=F_{0,\pm }(x,v).  \nonumber  \label{initial}
\end{eqnarray}
Here $F_{\pm }(t,x,v)\ge 0$ are the spatially periodic number
density functions for the ions (+) and electrons (-) respectively,
at time $t\ge 0$, position $x=(x_1,x_2,x_3)\in {\bf [}-\pi {\bf
,}\pi ]^3={\bf T}^3$, velocity $v=(v_1,v_2,v_3)\in {\bf R}^3$, and
$e_{\pm },$ $m_{\pm }$ the magnitude of their charges and masses.

 The collision between particles is given by
the standard Boltzmann collision operator $Q(G_1,G_2)$: Let $G_1(v),$ $G_2(v)$ be two number
density functions for two types of particles with masses $m_i$ and
diameters $\sigma _i$ ($i=1,2)$, then (p 83 and p 89 in \cite{Ch})
$Q(G_1,G_2)(v)$ is defined as
\begin{eqnarray}
&&\frac 14(\sigma _1+\sigma _2)^2\int_{{\bf R}^3\times
S^2}b(u-v, \omega) \{G_1(v^{\prime })G_2(u^{\prime
})-G_1(v)G_2(u)\}dud\omega
\label{hard} \\
&\equiv &Q_{{\rm {gain}}}(G_1,G_2)-Q_{{\rm {loss}}}(G_1,G_2).
\nonumber
\end{eqnarray}
For hard-sphere interaction, the collision kernel $b$ is given by $b(u-v, \omega)=|(u-v)\cdot w|$ (see \cite{Ch}) {and in this article, we assume the hard-sphere interaction, but the formal derivation will be valid for other general collision kernels.}
Here $\omega \in S^2$, and
\begin{equation}
v^{\prime }=v-\frac{2m_2}{m_1+m_2}[(v-u)\cdot \omega ]\omega ,\qquad
u^{\prime }=u+\frac{2m_1}{m_1+m_2}[(v-u)\cdot \omega ]\omega ,
\label{prime}
\end{equation}
which denote velocities after a collision of particles having velocities $%
v,u $ before the collision and vice versa. Notice that the elastic collision (\ref{prime})
implies
the conservation of momentum $m_1v+m_2u$ and energy $\frac 12m_1|v|^2+$ $%
\frac 12m_2|u|^2$ during the collision process. To clarify the
collisions between two types of particles, we use the following
notation:
\[
\begin{split}
Q^+\equiv Q(F_+,F_+)\,;\; Q^\pm\equiv Q(F_+,F_-)\,;\; Q^\mp\equiv
Q(F_-,F_+)\,;\;Q^-\equiv Q(F_-,F_-)\,.
\end{split}
\]
Note that $Q^+$ and $Q^-$ are the usual collision operators of one
species.

The self-consistent, spatially periodic electromagnetic field $%
[E(t,x),B(t,x)]$ in (\ref{VMB}) is coupled with $F(t,x,v)$
through the Maxwell system:
\begin{eqnarray*}
&&\mu_0\epsilon_0\partial _tE-\nabla \times B=-\mu_0 J=-\mu_0 \int_{{\bf R}%
^3}v\{e_{+}F_{+}-e_{-}F_{-}\}dv,\;\; \\
&&\partial _tB+\nabla \times E=0,\;\;\;\;\nabla \cdot B=0, \\
&&\nabla \cdot E =\frac{1}{\epsilon_0} \int_{{\bf R}^3}\{e_{+}F_{+}-e_{-}F_{-}%
\}dv,\; \\
&&E(0,x)=E_0(x),\;\;B(0,x)=B_0(x).
\end{eqnarray*}
Here  $\epsilon_0$ and $\mu_0$ are called the electric permittivity and the magnetic permeability of the plasma. And the speed of light $c$ is given by 
\[
c^2=\frac{1}{\mu_0\epsilon_0}. 
\]
It is well-known that for classical solutions to the
VMB system, the following conservation laws of
mass, total momentum (both kinetic and electromagnetic) and total
energy (both kinetic and electromagnetic) hold:
\begin{eqnarray*}
&&\frac d{dt}\int_{{\bf T}^3\times {\bf R}^3}m_{+}F_{+}(t)=0,\;\;\; \\
\; &&\frac d{dt}\int_{{\bf T}^3\times {\bf R}^3}m_{-}F_{-}(t)=0, \\
&&\frac d{dt}\left\{ \int_{{\bf T}^3\times {\bf R}%
^3}v(m_{+}F_{+}(t)+m_{-}F_{-}(t))+\frac 1{4\pi }\int_{{\bf
T}^3}E(t)\times
B(t)\right\} =0, \\
&&\frac d{dt}\left\{ \frac 12\int_{{\bf T}^3\times {\bf R}%
^3}|v|^2(m_{+}F_{+}(t)+m_{-}F_{-}(t))+{\frac 1{8\pi }}\int_{{\bf T}%
^3}|E(t)|^2+|B(t)|^2\right\} =0.
\end{eqnarray*}
Moreover, we also have the celebrated H-Theorem of
Boltzmann
\begin{equation}
\frac d{dt}\left\{ \int_{{\bf T}^3\times {\bf R}^3}(F_{+}(t)\ln
F_{+}(t)+F_{-}(t)\ln F_{-}(t))\right\} \le 0.  \label{entropy}
\end{equation}

Indeed, the above conservation laws and the entropy dissipation can be derived from the following well-known properties of the collision operators
$Q^+,\,Q^\pm,\,Q^\mp,\,Q^-$. We refer to \cite{ABT,Ch,Degond} for more details. 

\smallskip

(i) Mass conservation:
\[
\int Q^+ dv=0\,,\int Q^\pm dv=0\,,\int Q^\mp dv=0\,,\int Q^- dv=0\,.
\]

(ii) Momentum conservation:
\[
\int Q^+ m_+ v dv=0\,,\int (Q^\pm m_+v + Q^\mp m_-v) dv=0\,,\int Q^-
m_- vdv=0\,.
\]

(iii) Energy conservation:
\[
\int Q^+ m_+ |v|^2 dv=0\,,\int (Q^\pm m_+|v|^2 + Q^\mp m_-|v|^2)
dv=0\,,\int Q^- m_- |v|^2dv=0\,.
\]

(iv) Entropy inequalities:
\[
\begin{split}
\int Q^+(F_+,F_+) \ln F_+ dv\leq 0\,,\int Q^-(F_-,F_-) \ln F_-
dv\leq 0, \\
\int (Q^\pm(F_+,F_-) \ln F_+  + Q^\mp(F_-,F_+) \ln F_-) dv\leq 0. 
\end{split}
\]

As a consequence of the entropy inequalities, we can find
distribution functions (two Maxwellians) which cancel the collision
operators.

\smallskip

(v) Local Thermodynamical Equilibria:
\begin{equation}\label{LM}
M _{+}(t,x,v)={n_+}(\frac{m_{+}}{2\pi T}%
)^{3/2}e^{-m_{+}|v-\mathfrak{u}|^2/2 T},\;
M _{-}(t,x,v)={n_-}(\frac{%
m_{-}}{2\pi T})^{3/2}e^{-m_{-}|v-\mathfrak{u}|^2/2 T}
\end{equation}
so that
$$Q^+(M_+,M_+)+Q^\pm(M_+,M_-)=0\text{ and }Q^\mp(M_-,M_+)+Q^-(M_-,M_-)=0$$
where $n_+$ and $n_-$ are the ion and electron density and
$\mathfrak{u}$ and $T$ are the common mean velocity and temperature
and they may depend on $t$ and $x$. When the macroscopic variables do not depend on 
 $t$ and $x$,  we call them  global Maxwellians. When $\mathfrak{u} = 0 $, we define 
\begin{equation}\label{GM}
\mu _{+}(v)={n_+}(\frac{m_{+}}{2\pi  T_0}%
)^{3/2}e^{-m_{+}|v|^2/2 T_0},\;\;\;\mu _{-}(v)={n_-}(\frac{%
m_{-}}{2\pi  T_0})^{3/2}e^{-m_{-}|v|^2/2T_0}.
\end{equation} 
Note that  {when $e_+ {n_+} = e_-  {n_-} $}, 
  the  global Maxwellians  $\mu _{+}(v) $  and $\mu _{-}(v) $ 
  together with $E=B=0$ define a stationary 
 solutions to the VMB system, while the local Maxwellians are not necessarily 
 solutions.

We now introduce the linearized collision operators over the
velocity space for $Q^++Q^\pm$ and $Q^\mp+Q^-$ around $M_+$ and
$M_-$ in the following vector form:
\begin{equation}\label{L}
\begin{split}
&-\mathfrak{L}\binom{f_+}{f_-}\equiv\\ &\binom{M_+^{-1}\{Q^+(M_+f_+,M_+)+Q^+(M_+,
M_+f_+)+Q^\pm(M_+f_+,M_-)+Q^\pm(M_+, M_-f_-)\}}{M_-^{-1}\{Q^\mp(M_-f_-,M_+)+Q^\mp(M_-,
M_+f_+)+Q^-(M_-f_-,M_-)+Q^-(M_-, M_-f_-)\}}.  
\end{split}
\end{equation}
We define the following  inner
product in $v$: 
\[
\Big\langle \binom{f_+}{f_-},\,\binom{g_+}{g_-}\Big\rangle_M
=\int_{\mathbf{R}^3} (M_+f_+ g_++M_-f_-g_-)\, dv\,
\]
and denote by $L^2_M$ the associated Hilbert space. 
We will use $\langle \cdot,\cdot\rangle$ to denote the usual $L^2$ inner product without the weight. 
We summarize the properties of $\mathfrak{L}$ in the following lemma. The proof can be found  in \cite{ABT} 
(pp 635--638).

\smallskip

\begin{lemma}\label{lem}  We assume the hard-sphere interaction for the collision kernel. 

1.  $\mathfrak{L}$ is the sum of a diagonal operator $f\rightarrow \nu f$
$$ \nu f=\binom{\nu_+(v) f_+}{\nu_-(v)f_-} \;\text{ with } \;\nu_\pm(|v|)\sim 1+|v| $$
and a compact operator $\mathfrak{K}$.  The domain of $\mathfrak{L}$ is given by 
 $$D(\mathfrak{L})=\{f: \|  (1+|v|)^{\frac12} f \|_{L^2_M}<\infty\}.$$

2.  $\mathfrak{L}$ is self-adjoint in $L^2_M $:
$$\Big\langle \mathfrak{L}\binom{f_+}{f_-}, \binom{g_+}{g_-}\Big \rangle_M=
\Big\langle \binom{f_+}{f_-},\mathfrak{L}\binom{g_+}{g_-} \Big \rangle_M.$$
$\mathfrak{L}$ is non-negative.

3. The kernel of $\mathfrak{L}$ is a six-dimensional
linear space:
$$\ker \mathfrak{L}=\hbox{Span} \left\{\binom{1}{0},\,\binom{0}{1},\,
\binom{m_+ v_i}{m_-v_i},\,\binom{m_+|v|^2}{m_-|v|^2}\right\}.$$

4. Any function $f\in D(\mathfrak{L})$ can be written as $f= q_f+ w_f$ with 
$q_f\in \ker \mathfrak{L}$ and $w_f\in (\ker \mathfrak{L})^\perp$ and we have $\langle \mathfrak{L}f,f \rangle_M\geq \delta_0\|(1+|v|)^{\frac12} w_f \|^2$. 
\end{lemma}

 Note that this Fredholm operator $\mathfrak{L}$ can be inverted after checking that
the inhomogeneity is perpendicular to its six-dimensional null space
$$\left\{\binom{1}{0},\,\binom{0}{1},\, \binom{m_+
v_i}{m_-v_i},\,\binom{m_+|v|^2}{m_-|v|^2}\right\}\equiv
\{\phi_0,...,\phi_5\}.$$ 

\subsection{F-G formulation}\label{F-G}

In the MHD model, the medium is considered as a single fluid, in other words, one neglects the difference in motion of the electrons, various kinds of ions, and neutral particles. To derive the MHD equations from the VMB system \eqref{VMB}, it is convenient to consider the total mass density and the total charge density for $F^+$ and $F^-$: 
\begin{equation}
F\equiv m_{+}F_{+}+m_{-}F_{-}\;\text{ and } \;G\equiv e_+F_{+}-e_-F_{-}\,.\label{FG}
\end{equation}
Then the VMB system \eqref{VMB} can be rewritten as
\begin{equation}
\begin{split}\label{rVMB}
&\partial_t F+v\cdot\nabla_x F+(E+{v}\times B)\cdot
\nabla_vG 
 =m_+( Q^++Q^{\pm})+m_-(Q^\mp+Q^-), \\
&\partial_t G+v\cdot\nabla_x G +(E+{v}\times B)\cdot
\nabla_v \{ \frac{e_{+}e_{-}}{m_{+}m_{-}} F +
\frac{e_{+}m_{-}-e_{-}m_{+}}{m_{+}m_{-}} G\} \\&\quad\quad\quad\quad\quad\quad\quad\quad\quad\quad
\quad\quad\quad\quad\quad\;\;\;=e_+( Q^++Q^{\pm})-e_-(Q^\mp+Q^- ),
\end{split}
\end{equation}
and  
\begin{equation}
\begin{split}\label{rM}
\mu_0\epsilon_0\partial _tE-\nabla \times B=-\mu_0 \int_{{\bf R}%
^3}v\, Gdv,\;\; \partial _tB+\nabla \times E=0,\\
\nabla \cdot E=\frac{1}{\epsilon_0} \int_{{\bf R}^3} G dv, 
\;\;\;\;\nabla \cdot B=0.
\end{split}
\end{equation}
Notice that the Maxwell's equations are coupled with only $G$.  We will use both $F^+, F^-$ formulation \eqref{VMB} and $F,G$ formulation \eqref{rVMB} throughout the article. 
On one hand,  we remark that when $m_+=m_-(=m)$, the collision operators are of one-species and thus the right-hand-sides of \eqref{rVMB} can be significantly simplified: 
\begin{equation}
\begin{split}\label{fg}
\bullet\;&m_+( Q^++Q^{\pm})+m_-(Q^\mp+Q^-)\\
&= m\{Q(F_+,F_+) +Q(F_+,F_-)+Q(F_-,F_+)+Q(F_-,F_-) \}\\
&=m Q(F_++F_-,F_++F_-) \\
&= \frac{1}{m}Q(F,F)\\
\bullet\;& e_+( Q^++Q^{\pm})-e_-(Q^\mp+Q^- )\\
&= e_+\{Q(F_+,F_+) +Q(F_+,F_-)\}-e_-\{Q(F_-,F_+)+Q(F_-,F_-) \}\\
&=Q(e_+F_+-e_-F_-,F_++F_-)\\
&=\frac{1}{m}Q(G,F)
\end{split}
\end{equation}
In this case, the linearized operator \eqref{L} is orthogonally decomposed: the linearized operator of $Q(F,F)$ has the five-dimensional kernel $\{1,v,|v|^2\}$ and the linearized operator of $Q(G,F)$ has the one-dimensional kernel $\{1\}$, for instance see \cite{belm,j3} for the use of these linearized operators. This special setting will be frequently used. 
On the other hand, we remark that when $e_{+}m_{-} = e_{-}m_{+}  $, then the 
term $\frac{e_{+}m_{-}-e_{-}m_{+}}{m_{+}m_{-}} G  $ cancels in \eqref{rVMB}.


\section{Near local thermodynamical equilibria}

By using the collision invariants of the Boltzmann collision operators, if $F,G,E,B$ satisfy the VMB system \eqref{rVMB} and \eqref{rM}, we formally deduce the following local conservation laws: 
\begin{equation}
\begin{split}\label{local}
&\bullet\;\partial_t\Big(\int F dv\Big)+\nabla_x\!\cdot\!\Big(\int v\,F dv\Big) =0 \\
&\bullet\;\partial_t\Big(\int v\,F dv\Big) + \nabla_x\!\cdot\!\Big(\int v\otimes v\, Fdv\Big)=\Big(\int G dv\Big)E +\Big(\int v\,G dv\Big)\times B \\
&\bullet\;\partial_t \Big(\int \frac{|v|^2}{2}F dv\Big) + \nabla_x\!\cdot\!\Big(\int v\frac{|v|^2}{2}\, Fdv\Big) =E\!\cdot  \Big(\int v\,G dv \Big)
\\
&\bullet\;\partial_t \Big(\int G dv\Big) + \nabla_x\!\cdot\!\Big(\int v\,G dv\Big)=0
\end{split}
\end{equation}
which describe the conservation of the total mass, momentum, and energy and the conservation of the charge density. Of course, the equations are coupled with \eqref{rM}. However, the integrals inside $\nabla_x\cdot$ such as $\int v\otimes v\, Fdv, \int v\frac{|v|^2}{2}\, Fdv, \int v\,G dv$ are not functions of the other macroscopic quantities and they depend on the entire distributions $F$ and $G$ unless the appropriate forms for $F$ and $G$ are assumed. Thus, in general, the moment method does not give a closed form of a finite set of equations and this is often referred as the closure problem. In order to close the system of balance laws \eqref{local}, we look for the solutions very close to local thermodynamical equilibria  $F_+\sim M_+$ and $F_-\sim M_-$ such that in such asymptotic regime, the above integrals can be expressed in terms of the given macroscopic variables. This approximation can be realized by rescaling the system, namely by comparing the size of the system and the time scale of interest with the typical space and time scales of collisions. The task is not all trivial because the physical regimes are much richer than the pure Boltzmann equation and one needs to deal with multi-scales encoded in the VMB system. 

 In what follows, we introduce a few different hyperbolic scalings to the VMB system and study their hydrodynamic limits formally. Also, we discuss the Hall effect of the MHD equations.  We want to point out that we did not try to get all possible 
scalings. In particular, one can easily  derive simplified models to those 
we have here. 
 

\subsection{$1\frac12$-- fluid Euler-Maxwell system with the common velocity and temperature}

We start with the simplest possible hyperbolic scaling to \eqref{VMB}:
Let  $\tilde{t}=\varepsilon t$, $\tilde{x}=\varepsilon x$, $\tilde{F}_\pm=\frac{1}{\varepsilon^2}F_\pm$, $\tilde{b}=\varepsilon^2 b$,  $\tilde{E}= \frac{1}{\varepsilon} E$ and $\tilde{B}=\frac{1}{\varepsilon}B$. Dropping $\tilde{\text{ }} $, we obtain the following rescaled VMB system: 
\begin{equation}
\begin{split}\label{rvmbC}
\partial _tF_{+}^\varepsilon+v\cdot \nabla _xF_{+}^\varepsilon+\frac{e_{+}}{m_{+}}(E^\varepsilon+
v\times B^\varepsilon)\cdot \nabla _vF^\varepsilon_{+}=\frac{1}{\varepsilon}
\{Q^+(F^\varepsilon_{+},F^\varepsilon_{+})+Q^\pm(F^\varepsilon_{+},F^\varepsilon_{-})\}, \\
\partial _tF^\varepsilon_{-}+v\cdot \nabla _xF^\varepsilon_{-}-\frac{e_{-}}{m_{-}}(E^\varepsilon+
v\times B^\varepsilon)\cdot \nabla
_vF^\varepsilon_{-}=\frac{1}{\varepsilon}\{Q^\mp(F^\varepsilon_{-},F^\varepsilon_{+})+Q^-(F^\varepsilon_{-},F^\varepsilon_{-})\},\\
\mu_0\epsilon_0\partial _tE^\varepsilon-\nabla \times B^\varepsilon=-\mu_0 \int_{{\bf R}
^3}v\{e_{+}F^\varepsilon_{+}-e_{-}F^\varepsilon_{-}\}dv,\;\;
\partial _tB^\varepsilon+\nabla \times E^\varepsilon=0, \\
\nabla \cdot E^\varepsilon =\frac{1}{\epsilon_0}\int_{{\bf R}^3}\{e_{+}F^\varepsilon_{+}-e_{-}F^\varepsilon_{-}
\}dv,\;\;\nabla \cdot B^\varepsilon=0, 
\end{split}
\end{equation}
where $\varepsilon>0$ is a Knudsen number. 

Suppose $F_\pm^\varepsilon\rightarrow F_\pm$, $[E^\varepsilon,B^\varepsilon]\rightarrow [E,B]$ as $\varepsilon \rightarrow 0$. Then from the first two equations in \eqref{rvmbC}, we see that 
\[
Q^+(F^+,F^+)+Q^\pm(F^+,F^-)
=0 \,;\; Q^\mp(F^-,F^+)+Q^-(F^-,F^-)=0
\]
from which we deduce that
\[
F^+=M_+\;\text{ and }\; F^-=M_-
\]
local Maxwellians given in \eqref{LM}. There are six fluid variables $n_+$, $n_-$, $\mathfrak{u}$ and $T$ to be determined: the dynamics of  $n_+$, $n_-$, $\mathfrak{u}$ and $T$ can be derived by taking the moments of the first two equations. The system of local conservation laws 
\begin{equation}\label{solvability}
\begin{split}
\left\langle \binom{\partial_tF_+^\varepsilon+v\cdot \nabla_x
F_+^\varepsilon+\frac{e_+}{m_+}(E^\varepsilon+{v}\times B^\varepsilon)\cdot
\nabla_vF_+^\varepsilon}{\partial_tF_-^\varepsilon+v\cdot \nabla_x
F_-^\varepsilon-\frac{e_-}{m_-}(E^\varepsilon+{v}\times B^\varepsilon)\cdot \nabla_v
F_-^\varepsilon}, \phi_i\right\rangle =0\,,\;i=0,...,5
\end{split}
\end{equation}
is, in general, not closed. If each moment converges, we can pass to the limit of $\varepsilon\rightarrow 0$. 
To write the equations for $n_+$, $n_-$, $\mathfrak{u}$ and $T$, we first compute the first few moments of local Maxwellians $M_\pm$:  
\[
\begin{split}
&\int M_\pm dv= {n_\pm}, \quad \;\int vM_\pm dv= {n_\pm}\mathfrak{u},\quad \;
\int |v|^2 M_\pm dv= {n_\pm}|\mathfrak{u}|^2 + 3\frac{n_\pm T}{m_\pm }, \\
& \int v\otimes v M_\pm dv = {n_\pm}\mathfrak{u}\otimes\mathfrak{u} + \frac{n_\pm T}{m_\pm } \delta^i_j , \quad \int |v|^2 v M_\pm dv= {n_\pm}  |\mathfrak{u}|^2\mathfrak{u} + 5   \frac{n_\pm  T}{m_\pm } \mathfrak{u}. 
\end{split}
\]
Thus from \eqref{solvability} where $F^\varepsilon_\pm$ are placed by $F_\pm$
\begin{equation}
\begin{split}\label{EM0}
&\bullet\;\partial_tn_+ + \nabla_x\!\cdot\!(n_+\mathfrak{u})=0\\
&\bullet\;\partial_tn_- + \nabla_x\!\cdot\!(n_-\mathfrak{u})=0\\
&\bullet\;\partial_t\Big[({m_+n_+}+{m_-n_-})\mathfrak{u}\Big]+\nabla_x\!\cdot\!\Big[({m_+n_+}+{m_-n_-}) \mathfrak{u}\otimes\mathfrak{u}\Big]+ \nabla_x\Big[({n_+}+{n_-}) T)\Big]\\
&\quad= (n_+e_+-n_-e_-) E+{(n_+e_+-n_-e_-)\mathfrak{u}\times B}\\
&\bullet\;\partial_t\Big[\frac12{({m_+n_+}+{m_-n_-}) |\mathfrak{u}|^2}+ \frac{3}{2} {({n_+}+{n_-})}T\Big] \\
&\quad+\nabla_x\!\cdot\!\Big[ \frac{1}{2}{({m_+n_+}+{m_-n_-})|\mathfrak{u}|^2\mathfrak{u}} + \frac{5}{2}  ({n_+}+{n_-}) T \mathfrak{u}  \Big]=E\!\cdot\!\Big[ ( n_+e_+-n_-e_- )\mathfrak{u} \Big]\\
\end{split}
\end{equation}
with 
\begin{equation}\label{EM00}
\begin{split}
\mu_0\epsilon_0\partial_tE-\nabla\times B=-\mu_0(e_+n_+-e_-n_-)\mathfrak{u};\;
\partial_tB+\nabla\times E=0;\\ 
\nabla \cdot E=\frac{1}{\epsilon_0}(e_+n_+-e_-n_-);\;\nabla\cdot B=0.
\end{split}
\end{equation}
By introducing the total density $\rho$ and the charge density $\sigma$: 
\[
{m_+n_+}+{m_-n_-} =: \rho \text{ and } e_+n_+-e_-n_- =:\sigma
\]
the above system \eqref{EM0} and \eqref{EM00}  can be equivalently written as 
\begin{equation}
\begin{split}\label{EM1}
&\bullet\partial_t\rho + \nabla_x\!\cdot\!(\rho\mathfrak{u})=0\\
&\bullet\partial_t(\rho\mathfrak{u})+\nabla_x\!\cdot\!(\rho \mathfrak{u}\otimes\mathfrak{u})+ \nabla_x\Big[({n_+}+{n_-}) T\Big]= \sigma E+\sigma\mathfrak{u}\times B\\
&\bullet\partial_t\Big[\frac{\rho |\mathfrak{u}|^2}{2} + \frac{3}{2} {({n_+}+{n_-})} T\Big]+\nabla_x\!\cdot\!\Big[ \frac{\rho|\mathfrak{u}|^2\mathfrak{u}}{2} + \frac{5}{2}  ({n_+}+{n_-}) T \mathfrak{u}  \Big] =\sigma E\!\cdot\! \mathfrak{u}\\
&\bullet\partial_t\sigma + \nabla_x\!\cdot\!(\sigma\mathfrak{u})=0
\end{split}
\end{equation}
where 
\begin{equation}\label{EM2}
\mu_0\epsilon_0\partial_tE-\nabla\times B=-\mu_0\sigma\mathfrak{u}\,;\;
\partial_tB+\nabla\times E=0\,;\; 
\nabla \cdot E=\frac{1}{\epsilon_0}\sigma\,;\;\nabla\cdot B=0
\end{equation}
and 
\begin{equation}\label{pr} 
{n_\pm} = \frac{e_\mp \rho \pm  m_\mp \sigma }{e_-m_++e_+m_- } 
\ \hbox{and} \ 
{n_+}+{n_-}= \frac{(e_++e_-)\rho+(m_--m_+)\sigma}{e_-m_++e_+m_-}\,.
\end{equation}
The first three equations in \eqref{EM1} represent the balance laws of total mass, total momentum, and total energy and the last equation is the conservation of the charge density. In this regime, the current density $J$ is given by $\sigma\mathfrak{u}$. Hence, by letting 
\[
p=({n_+}+{n_-}) T \;\text{ (ideal gas law)}, 
\]
 the solvability condition \eqref{solvability} yields Euler-Maxwell type equations \eqref{EM1} and \eqref{EM2} having the same mean velocity and temperature fields which are therefore not exactly of two-fluid.   From \eqref{pr}, note that for a special case of $m_+=m_-$, the pressure is induced by the total density but in general, it depends on the charge density as well. 

We have proved the following theorem. 

\begin{theorem} Let $F_\pm^\varepsilon$, $[E^\varepsilon,B^\varepsilon]$ be the classical solutions to \eqref{rvmbC}. Assume that  $F_\pm^\varepsilon\rightarrow F_\pm$, $[E^\varepsilon,B^\varepsilon]\rightarrow [E,B]$ strongly as $\varepsilon \rightarrow 0$. Then the limit $F_\pm$ is a Maxwellian given by \eqref{LM} and moreover, $n_\pm,\mathfrak{u},T$, $E,B$ solve the $1\frac12$-fluid Euler-Maxwell system 
\eqref{EM0} and \eqref{EM00}. 
\end{theorem}

For a rigorous proof of this result, we can argue as in 
Caflisch \cite{C}. Indeed, the derivation given above can also 
be used to provide a Hilbert expansion as in  \cite{C}. 
This is then used to provide an error estimate between 
the solution to the  $1\frac12$-fluid Euler-Maxwell system and 
the solutions to the VMB system \eqref{rvmbC} in the limit 
when $\varepsilon$ goes to zero. 


\subsection{Compressible resistive MHD}

We have seen that it is convenient to write the equations in terms of the total density and the charge density in the previous section. In this subsection, we start with the $F-G$ formulation introduced in Section \ref{F-G}. We set $m_+=m_-=e_+=e_-=1$ for simplicity of the presentation. Then by using \eqref{fg}, the VMB system \eqref{rVMB} reads as follows: 
\begin{equation}
\begin{split}\label{VMBquasi}
&\partial_t F+v\cdot\nabla_x F+(E+{v}\times B)\cdot
\nabla_vG 
 =Q(F,F), \\
&\partial_t G+v\cdot\nabla_x G +(E+{v}\times B)\cdot
\nabla_v F =Q(G,F).
\end{split}
\end{equation}
In order to see the ideal MHD type equations with  Ohm's law, we consider the 
 quasi-neutral regime so that $G\sim 0$. We next introduce another scaling to  \eqref{VMBquasi} as follows: Let  $\tilde{t}=\varepsilon t$, $\tilde{x}=\varepsilon x$,  $\tilde{E}= \frac{1}{\sqrt{\varepsilon}} E$, $\tilde{B}=\frac{1}{\sqrt{\varepsilon}}B$, $\tilde{G}=\frac{1}{{\varepsilon}^{3/2}}G$, $\tilde{F}=\frac{1}{{\varepsilon}}F$, and $\tilde{b}=\varepsilon b$ where $E=O(\sqrt{\varepsilon}),$ $B=O(\sqrt{\varepsilon})$, $G=O({\varepsilon}^{3/2})$ and $F=O(\varepsilon)$. Dropping $\tilde{\text{ }}$, we obtain the following rescaled VMB system: 
\begin{equation}
\begin{split}\label{rVMBquasi3}
&\partial_t F^\varepsilon+v\cdot\nabla_x F^\varepsilon+(E^\varepsilon+{v\times {B^\varepsilon}})\cdot
\nabla_vG ^\varepsilon
 =\frac{1}{\varepsilon}Q(F^\varepsilon,F^\varepsilon), \\
&\partial_t G^\varepsilon+v\cdot\nabla_x G^\varepsilon + (\frac{E^\varepsilon}{\varepsilon}+\frac{v\times {B^\varepsilon}}{\varepsilon})\cdot
\nabla_v F^\varepsilon =\frac{1}{\varepsilon}Q(G^\varepsilon,F^\varepsilon),
\end{split}
\end{equation}
coupled with the Maxwell's equations 
\begin{equation}\label{rM3}
\begin{split}
\mu_0\epsilon_0 \partial_t E^\varepsilon- \nabla\times B^{\varepsilon}=-\mu_0\int vG^\varepsilon dv,  \;\;\partial_tB^\varepsilon+\nabla\times E^\varepsilon =0,\\
 \nabla\cdot E^\varepsilon= \frac{1}{\epsilon_0}\int G^\varepsilon dv,\;\;\nabla\cdot B^\varepsilon=0. 
\end{split}
\end{equation}
Suppose that there exists  a  classical solution $(F^\varepsilon, G^\varepsilon, E^\varepsilon, B^\varepsilon)$ to \eqref{rVMBquasi3} and \eqref{rM3} and that there exist $F,G,E,B$ so that $F^\varepsilon\rightarrow F, G^\varepsilon\rightarrow G, E^\varepsilon\rightarrow E, B^\varepsilon\rightarrow B$ as $\varepsilon\rightarrow 0$. 
Then from \eqref{rVMBquasi3}, as $\varepsilon\rightarrow 0$ we deduce that
\begin{equation}\label{Flm}
Q(F,F)=0\;\;\Longrightarrow \;\; F=M=\frac{\rho }{(2\pi T)^{3/2}}  e^{-|v-\mathfrak{u}|^2/2T} 
\end{equation}
also 
\begin{equation}\label{Flm2}
(E+v\times B )\cdot\nabla_v F = Q(G,F). 
\end{equation}
This case, it is already different from the previous ones in that $G$ is not 
determined by a Maxwellian. Of course here $G$ is a higher order fluctuation. 
 Define 
\[
\mathcal{L} g=-\frac{1}{M} Q(Mg,M)
\]
where $M$ is the Maxwellian given in \eqref{Flm}. Then from Lemma \ref{lem} we deduce that $\ker \mathcal{L} =\text{span} \{1\}$. The solvability condition for $G$ is automatically satisfied and thus by inverting the linearized operator $\mathcal{L}$ in \eqref{Flm2} we derive that 
\begin{equation}
\begin{split}\label{ohm0}
G = -M\big[\mathcal{L}^{-1}\{\frac{1}{M}(E+v\times B )\cdot\nabla_v M\}+\psi_0\big], \text{ where }\mathcal{L}\psi_0=0 . 
\end{split} 
\end{equation}
Let us define the first two moments of $G$ by $\sigma$ (the charge density) 
  and $J$ (the current): 
\begin{equation}\label{J0}
\sigma:=\int G dv\;\text{ and }\;J:=\int v G dv.  
\end{equation}
Then from \eqref{ohm0} we can write the current $J$ as 
\[
J=\int vG dv=\langle v, G \rangle
 = -\langle v,M \mathcal{L}^{-1}\{\frac{1}{M}(E+v\times B )\cdot\nabla_v M\} \rangle 
\]
and since 
\[
\frac{1}{M}(E+v\times B )\cdot\nabla_v M = - (E+v\times B )\cdot \frac{(v-\mathfrak{u})}{T} =- \frac{1}{T} ( E+\mathfrak{u}\times B )\cdot (v-\mathfrak{u}) 
\]
and by using the fact that $\mathcal{L}^{-1}$ acts only in $v$, we obtain 
\[
\begin{split}
J&= \langle v,M \mathcal{L}^{-1}\{\frac{1}{T} ( E+\mathfrak{u}\times B )\cdot (v-\mathfrak{u}) \} \rangle = 
 \langle vM, \frac{1}{T}  ( E+\mathfrak{u}\times B )\cdot  \mathcal{L}^{-1}(v-\mathfrak{u}) \rangle \\
 &=  \langle ( v-\mathfrak{u})M,  \frac{1}{T}  ( E+\mathfrak{u}\times B )\cdot  \mathcal{L}^{-1}(v-\mathfrak{u}) \rangle. 
\end{split}
\]
Since  $\int   (v^i-\mathfrak{u}^i) \mathcal{L}^{-1}(v^j-\mathfrak{u}^j) M dv =0 $ for $i\neq j$, the last expression for $J$ can be reduced to the following: 
\begin{equation}
\begin{split}\label{ohm1}
J&= \frac{1}{\eta} (E+\mathfrak{u}\times B)
\end{split} 
\end{equation}
where 
\[
\frac{1}{\eta}=\frac{1}{3T}\int (v-\mathfrak{u})\cdot  \mathcal{L}^{-1}(v-\mathfrak{u}) M dv >0 . 
\]
Here, the number $3$ represents the space dimension 3. The relation \eqref{ohm1} is called the Ohm's law.  $\eta$ is called  \textit{resistivity} and $1/\eta$ is the conductivity. Next we take the moments of $F$ equation in \eqref{rVMBquasi3}, then by collision invariants of $Q$, we get 
\begin{equation}
\begin{split}\label{FG_0}
&\langle\partial_t F^\varepsilon+v\cdot\nabla_x F^\varepsilon+(E^\varepsilon+v\times B^\varepsilon)\cdot
\nabla_vG^\varepsilon\,, \{1,v,\frac{|v|^2}{2}\}\rangle
 =0
\end{split}
\end{equation}
and 
\begin{equation}\label{M_0}
\begin{split}
\mu_0\epsilon_0\partial_tE^\varepsilon- \nabla\times B^\varepsilon=- \mu_0\int vG^\varepsilon dv,  \;\; \partial_tB^\varepsilon+\nabla\times E^\varepsilon =0,\\
\nabla\cdot E^\varepsilon = \frac{1}{\epsilon_0}\int G^\varepsilon dv,\;\;\nabla\cdot B^\varepsilon=0. 
\end{split}
\end{equation}
We pass to the limit and by using the fact that $F$ is a Maxwellain for which we can explicitly compute the moments and by using  the notation introduced in \eqref{J0}, 
 we obtain the following resistive MHD type equations for $\rho, \mathfrak{u}, T, \sigma, J, E,B$: 
\begin{equation}
\begin{split}\label{mhd}
&\bullet\;\partial_t\rho + \nabla_x\!\cdot\!(\rho\mathfrak{u})=0\\
&\bullet\;\partial_t(\rho\mathfrak{u})+\nabla_x\!\cdot\!(\rho \mathfrak{u}\otimes\mathfrak{u})+ \nabla_x(\rho T)=\sigma E+ {J\times B}\\
&\bullet\;\partial_t\Big(\frac{\rho |\mathfrak{u}|^2}{2} + \frac{3}{2}\rho T\Big)+\nabla_x\!\cdot\!\Big( \frac{\rho|\mathfrak{u}|^2\mathfrak{u}}{2} + \frac{5}{2} \rho T \mathfrak{u}  \Big)=E\!\cdot \! J\\
&\bullet\;\partial_t\sigma + \nabla_x\!\cdot\!J=0\\
&\bullet\; J=  \frac{1}{\eta} (E+\mathfrak{u}\times B)\\
&\bullet\;  \mu_0\epsilon_0 \partial_tE-\nabla\times B=-\mu_0J ,  \;\; \partial_tB+\nabla\times E =0,\;\; \nabla\cdot B=0, \; \nabla \cdot E=\frac{\sigma}{\epsilon_0}.  
\end{split}
\end{equation}
We have proved the following theorem. 

\begin{theorem} Let $[F^\varepsilon,G^\varepsilon]$, $[E^\varepsilon,B^\varepsilon]$ be the classical solutions to \eqref{rVMBquasi3} and \eqref{rM3}. Assume that  $[F^\varepsilon,G^\varepsilon]\rightarrow [F,G]$, $[E^\varepsilon,B^\varepsilon]\rightarrow[E,B]$ strongly as $\varepsilon \rightarrow 0$. Then the limit $F$ is a Maxwellian given by \eqref{Flm} and moreover, $\rho,\mathfrak{u},T$, $\sigma, J$, $E,B$ solve the electromagneto-hydrodynamic system 
\eqref{mhd}. 
\end{theorem}

The electromagnetic field in the system \eqref{mhd} is still governed by the full Maxwell's equations; if we  consider the regime where  the electric permittivity $\epsilon_0$ is sufficiently small, namely if we consider the quasi-neutral limit, 
 then we would get a new system in which the magnetic field becomes dominant. 
We may set $\epsilon_0\sim \varepsilon$ by fixing $\mu_0$ in the Maxwell's equations \eqref{rM3} to obtain different approximations to the dynamics of the electromagnetic field. Namely, if we consider $$\mu_0\varepsilon \partial_t E^\varepsilon- \nabla\times B^{\varepsilon}=-\mu_0\int vG^\varepsilon dv\; \text{ and  }\;\nabla\cdot E^\varepsilon=\frac{1}{\varepsilon}\int G^\varepsilon dv,$$ following the same procedure, we would get $\nabla\times B=\mu_0J$ and $\sigma=0$. And \eqref{mhd} reduces to  
\begin{equation}
\begin{split}\label{mhd1}
&\bullet\;\partial_t\rho + \nabla_x\!\cdot\!(\rho\mathfrak{u})=0\\
&\bullet\;\partial_t(\rho\mathfrak{u})+\nabla_x\!\cdot\!(\rho \mathfrak{u}\otimes\mathfrak{u})+ \nabla_x(\rho T)= {J\times B}\\
&\bullet\;\partial_t\Big(\frac{\rho |\mathfrak{u}|^2}{2} + \frac{3}{2}\rho T\Big)+\nabla_x\!\cdot\!\Big( \frac{\rho|\mathfrak{u}|^2\mathfrak{u}}{2} + \frac{5}{2} \rho T \mathfrak{u}  \Big)=E\!\cdot\! J \\
&\bullet\; J=  \frac{1}{\eta} (E+\mathfrak{u}\times B)\\
&\bullet\;  \nabla\times B=\mu_0J ,  \;\; \partial_tB+\nabla\times E=0, \;\; \nabla\cdot B=0
\end{split}
\end{equation}
which are the well-known compressible resistive MHD equations.  

\begin{theorem} Let $[F^\varepsilon,G^\varepsilon]$, $[E^\varepsilon,B^\varepsilon]$ be the classical solutions to \eqref{rVMBquasi3} and \eqref{rM3} where $\epsilon_0$ is taken to be $\varepsilon$. Assume that  $[F^\varepsilon,G^\varepsilon]\rightarrow [F,G]$, $[E^\varepsilon,B^\varepsilon]\rightarrow[E,B]$ strongly as $\varepsilon \rightarrow 0$. Then the limit $F$ is a Maxwellian given by \eqref{Flm} and moreover, $\rho,\mathfrak{u},T, J$ and $E,B$ solve the compressible resistive MHD equations  \eqref{mhd1}. 
\end{theorem}

We remark that the ideal MHD equations are valid when the resistivity $\eta$ is negligible.  
 
Also, we point out that the simplification of setting the physical constants to be one allows us to treat the collision operators decoupled as in \eqref{VMBquasi}. The 
 General case, in particular when  $m_+\neq m_-$, can be treated similarly. In such case, 
we have to start with \eqref{VMB} or \eqref{rVMB}  and need to treat the collision operators as matrix form as in \eqref{L} and use that it has a six dimensional kernel. 
We do not detail this here.


\subsection{When $e_-m_+\neq e_+ m_-$ and Hall effect}

 The Hall effect is known as two-fluid effect. In our setting, the Hall effect can be captured when $e_-m_+\neq e_+ m_-$ at the \textit{higher order}. For instance, we take the same scaling as done in  \eqref{rVMBquasi3} to \eqref{rVMB},  
\begin{equation}
\begin{split}\label{rVMBh}
&\partial_t F^\varepsilon+v\cdot\nabla_x F^\varepsilon+(E^\varepsilon+{v}\times B^\varepsilon)\cdot
\nabla_vG^\varepsilon 
 =\frac{1}{\varepsilon} \{m_+( Q^++Q^{\pm})+m_-(Q^\mp+Q^-)\}, \\
&\partial_t G^\varepsilon+v\cdot\nabla_x G^\varepsilon +(\frac{E^\varepsilon}{\varepsilon}+ \frac{{v}\times B^\varepsilon}{\varepsilon})\cdot
\nabla_v \{ \frac{e_{+}e_{-}}{m_{+}m_{-}} F^\varepsilon +\varepsilon^{\frac12} 
\underbrace{\frac{e_{+}m_{-}-e_{-}m_{+}}{m_{+}m_{-}} G^\varepsilon}_{(\ast)}\} \\&\quad\quad\quad\quad\quad\quad\quad\quad\quad\quad
\quad\quad\quad\quad\quad\;\;\;=\frac{1}{\varepsilon} [e_+( Q^++Q^{\pm})-e_-(Q^\mp+Q^- )]. 
\end{split}
\end{equation}
The main difference is the presence of the term $(\ast)$, which will affect the Ohm's law. 

Following the same spirit of the Chapman-Enskog expansion, if we also  keep  the 
O($\sqrt{\varepsilon}$) terms in the second equation of 
\eqref{rVMBh},  instead of \eqref{ohm1}, 
 we would get the following Ohm's law
\begin{equation}\label{ohm2}
\eta J = \frac{e_{+}e_{-}}{m_{+}m_{-}} (  E + \mathfrak{u}\times B ) +\frac{\varepsilon^{\frac12} }{\rho}(\frac{e_{+}m_{-}-e_{-}m_{+}}{m_{+}m_{-}} ) (\sigma E+J\times B)
\end{equation}
where $J\times B$ is often referred to as  the Hall  effect 
term. We remark that this Hall effect is not seen by the moment method that we have used for the previous theorems but rather this can be captured by the higher order Chapman-Enskog type expansions.  It would be interesting to justify this higher order effect in a rigorous framework.

\begin{remark} We have chosen our scalings so far based on the appropriate time, space and forcing with fixed masses $m_+$ and $m_-$. When $m_-\ll m_+$ and $\frac{m_-}{m_+}$ is comparable to Knudsen number, one can also parametrize $\frac{m_-}{m_+}$. By doing so, for example, we can get the rescaled system \eqref{rVMBquasi3} or \eqref{rVMBh} without rescaling forcing terms. In such regime, one may want to study the appropriate collision operators as suggested in \cite{Degond}. 
\end{remark}


\section{The incompressible regime}

As we have seen in the previous section, the derivation of hydrodynamic equations from the Boltzmann equation
is closely related to finding approximate solutions of the Boltzmann
equation, because the fluid variables are defined and change in
space and time scales that are very large when measured in the units
of the mean free path and mean free time between collisions. In the present section, we'd like to discuss  the longer time scaling. We will mainly work on the simplified case where the constants are taken to be unity and present  general cases at the end when we discuss the Hall effect.

\subsection{Incompressible Navier-Stokes-Fourier-Maxwell equations} 

Recall the system \eqref{VMBquasi} with the Maxwell's equations: 
\begin{equation}
\begin{split}
&\partial _tF+v\cdot \nabla _xF+(E+v\times B)\cdot \nabla
_v G=Q(F, F),\\
&\partial _tG+v\cdot \nabla _xG+(E+v\times B)\cdot \nabla
_v F=Q(G,F),\\
&\mu_0\epsilon_0\partial _tE-\nabla \times B=-\mu_0\int_{{\mathbb{ R}}
^3}v\; G\;dv,\;\;\nabla \cdot B=0,\\
&\partial _tB+\nabla \times E=0,\;\;\nabla \cdot E=\frac{1}{\epsilon_0}\int_{{\mathbb{
R}}^3}G\;dv, \label{vmb}
\end{split}
\end{equation}
where $F=F_++F_-$ is the density for the
whole particles and $G=F_+-F_-$ represents the disparity between two
species (the charge density).

We will discuss how to obtain the incompressible Navier-Stokes-Fourier-Maxwell system. To keep the hyperbolic structure of Maxwell's equations, we take $\mu_0$ and $\epsilon_0$ as $O(1)$ and let $\tilde{t}=\varepsilon t$, $\tilde{x}=\varepsilon x$,  $\tilde{v}=\varepsilon v$ (longer time scale), $\tilde{E}= {{\varepsilon}} E$, $\tilde{B}=\varepsilon B$, and $\tilde{F}={{\varepsilon}} F$. Then the system \eqref{vmb} can be rescaled as follows: by dropping $\tilde{ \text{ }}$, 
\begin{equation}
\begin{split}
&\varepsilon\partial _tF^{\varepsilon}+v\cdot \nabla
_xF^{\varepsilon}+(\varepsilon E^{\varepsilon}+{v\times
B^{\varepsilon}})\cdot \nabla
_v G^{\varepsilon}=\frac{1}{\varepsilon}Q(F^{\varepsilon}, F^{\varepsilon}),\\
&\varepsilon\partial _tG^{\varepsilon}+v\cdot \nabla
_xG^{\varepsilon}+(\frac{E^{\varepsilon}}{\varepsilon}+\frac{v\times
B^{\varepsilon}}{\varepsilon^2})\cdot \nabla
_v F^{\varepsilon}=\frac{1}{\varepsilon}Q(G^{\varepsilon},F^{\varepsilon}),\\
&\mu_0\epsilon_0\partial _tE^{\varepsilon}-\nabla \times
B^{\varepsilon}=-{\mu_0} \int_{{\mathbb{ R}}
^3}v \;G^{\varepsilon}\;dv,\;\;\nabla \cdot B^{\varepsilon}=0,\\
&\partial _tB^{\varepsilon}+\nabla \times
E^{\varepsilon}=0,\;\;\nabla \cdot
E^{\varepsilon}=\frac{1}{\epsilon_0}\int_{{\mathbb{
R}}^3}G^{\varepsilon}\;dv.\label{rvmbM}
\end{split}
\end{equation}
We now consider the following perturbation of $F^\varepsilon,G^\varepsilon, E^\varepsilon, B^\varepsilon$ around $[\mu,0,0,0]$: 
\[
F^\varepsilon =\mu (1+\varepsilon f_1^\varepsilon), \;G^\varepsilon=\varepsilon \mu g_1^\varepsilon, \;
E^\varepsilon =\varepsilon E_1^\varepsilon,\; B^\varepsilon =\varepsilon B_1^\varepsilon
\]
where $\mu$ is a global Maxwellian normalized as $$\frac{1}{(2\pi)^{3/2}} e^{-|v|^2/2}.$$
Suppose that there exist the classical solutions $f_1^\varepsilon,\,g_1^\varepsilon,\,E_1^\varepsilon,\, B_1^\varepsilon$ to \eqref{rvmbM} and further assume that as $\varepsilon\rightarrow 0$, $f_1^\varepsilon\rightarrow f$, $g_1^\varepsilon\rightarrow g$, 
$E_1^\varepsilon\rightarrow E$, $B_1^\varepsilon\rightarrow B$. Under this convergence assumption, we will derive the dynamics of the limiting variables $f,g,E,B$. We first rewrite the system \eqref{rvmbM} in terms of the perturbed variables: 
\begin{equation}
\begin{split}
&\varepsilon\partial _tf_1^{\varepsilon}+v\cdot \nabla
_xf_1^{\varepsilon}+\varepsilon (\varepsilon E_1^{\varepsilon}+{v\times
B_1^{\varepsilon}})\cdot \frac{1}{\mu}\nabla
_v (\mu g_1^{\varepsilon}) \\
&\quad\quad\quad\quad\quad\quad\quad\quad\quad=\frac{1}{\varepsilon \mu}\{ Q(\mu f_1^\varepsilon, \mu) +Q(\mu,\mu f_1^\varepsilon ) +\varepsilon Q(\mu f_1^\varepsilon,\mu f_1^\varepsilon)\}\\
&\varepsilon\partial _tg_1^{\varepsilon}+v\cdot \nabla
_xg_1^{\varepsilon}+(\frac{E^{\varepsilon}}{\varepsilon}+\frac{v\times
B^{\varepsilon}}{\varepsilon^2})\cdot  \frac{1}{\mu}\nabla
_v(\mu+\varepsilon\mu f_1^{\varepsilon}) \\
&\quad\quad\quad\quad\quad\quad\quad\quad\quad=\frac{1}{\varepsilon \mu} \{Q(\mu g_1^\varepsilon , \mu) +\varepsilon Q(\mu g_1^\varepsilon, \mu f_1^{\varepsilon})  \},\\
&\mu_0\epsilon_0\partial _tE_1^{\varepsilon}-\nabla \times
B_1^{\varepsilon}=-{\mu_0} \int_{{\mathbb{ R}}
^3}v \mu g_1^{\varepsilon}\;dv,\;\;\nabla \cdot B_1^{\varepsilon}=0,\\
&\partial _tB_1^{\varepsilon}+\nabla \times
E_1^{\varepsilon}=0,\;\;\nabla \cdot
E_1^{\varepsilon}=\frac{1}{\epsilon_0}\int_{{\mathbb{
R}}^3}\mu g_1^{\varepsilon}\;dv.\label{rvmbMp}
\end{split}
\end{equation}
We recall the linearized collision operators $L$ and $\mathcal{L}$ (see for instance \cite{belm,j3}): 
\[
\begin{split}
Lf= -\frac{1}{\mu}\{Q(\mu f,\mu)+Q(\mu,\mu f)\} \;\text{ and }\; \mathcal{L}g=-\frac{1}{\mu} Q(\mu g,\mu). 
\end{split}
\]
Note that $[L,\mathcal{L}]$ is equivalent to the linearized collision operator introduced in \eqref{L}. We also recall that $\ker L=\text{span}\{1,v,|v|^2\}$ and $\ker \mathcal{L}=\text{span}\{1\}$. 

By taking $\varepsilon \rightarrow 0$ in  the first two equations in \eqref{rvmbMp} and using the convergence assumption, we first obtain 
\begin{equation}\label{3.4}
0=-Lf 
\end{equation}
and 
\begin{equation}\label{3.5}
-E\cdot v+\frac{1}{\mu}(v\times B)\cdot \nabla_v(\mu f) = - \mathcal{L}g\,.
\end{equation}
Since $f$ is in the kernel of $L$ by \eqref{3.4}, $f$ can be written as  $$f=\rho+v\cdot u+(\frac{|v|^2}{2}-\frac32)\theta\,.$$ 
We also introduce 
\[
\sigma := \int g \mu dv \;\text{ and } \;J:= \int g v\mu dv. 
\]
From \eqref{3.5}, we can write the current $J$ as 
\[
J=\int g v\mu dv = \int \mathcal{L}^{-1}(E\cdot v-\frac{1}{\mu}(v\times B)\cdot\nabla_v(\mu f)) v\mu dv. 
\]
Since 
\[
-\frac{1}{\mu}(v\times B)\cdot\nabla_v(\mu f) = - v\times B \cdot \nabla_v f = -v\times B \cdot u = u\times B \cdot v 
\]
and by using the fact that $\mathcal{L}^{-1}$ acts only in $v$, we obtain 
\[
J=\int (\mathcal{L}^{-1} v) \cdot (E+u\times B)\, v \mu dv . 
\]
Since $\int ( \mathcal{L}^{-1} v^i)  v^j \mu dv=0 $ for $i\neq j$, we can derive the following Ohm's law 
\begin{equation}\label{3.6}
\begin{split}
J= \frac{1}{\eta} (E+u\times B)
\end{split}
\end{equation}
where 
\[
\frac{1}{\eta} =\frac{1}{3}\int  (\mathcal{L}^{-1} v ) \cdot   v \mu dv . 
\]
The derivation of the dynamics of $\rho,u,\theta$ is very similar to the derivation of the incompressible Navier-Stokes-Fourier system  (for instance see \cite{BGL1}) except the presence of $E,B$ and $g$. 
The local conservation laws of mass and momentum of the first equation in \eqref{rvmbMp} are written as
\[
\begin{split}
&\varepsilon \langle \partial_t f^\varepsilon_1, \mu \rangle + \nabla_x\!\cdot \!\langle  f^\varepsilon_1, v\mu \rangle=0 \\ 
&\varepsilon \langle \partial_t f^\varepsilon_1, v\mu \rangle + \nabla_x\!\cdot\! \langle  f^\varepsilon_1,v\otimes v\mu \rangle + \varepsilon^2 \langle E_1^\varepsilon \cdot\nabla_v(\mu g_1^\varepsilon), v  \rangle + \varepsilon\langle  (v\times B_1^\varepsilon) \cdot\nabla_v(\mu g_1^\varepsilon),v \rangle=0 
\end{split}
\]
where we have used  $\int  (v\times
B_1)\cdot\nabla_v({\mu}g_1) dv =0$. 
By letting $\varepsilon$ go to zero, we see that 
\[
\nabla_x\!\cdot\! \langle  f, v\mu \rangle=0  \;   \text{ and } \;  \nabla_x\!\cdot\! \langle  f,v\otimes v\mu \rangle=0
\]
which yield the divergence-free condition and the Boussinesq relation
\[
\nabla_x\!\cdot\! u=0\; \text{ and }\; \nabla_x(\rho+\theta)=0. 
\]
The limiting momentum equation is obtained from 
\[
\underbrace{ \langle \partial_t f^\varepsilon_1, v\mu \rangle}_{(a)} +\underbrace{\frac{1}{\varepsilon} \nabla_x\!\cdot\! \langle  f^\varepsilon_1,v\otimes v\mu \rangle}_{(b)} +\underbrace{ \varepsilon \langle E_1^\varepsilon \cdot\nabla_v(\mu g_1^\varepsilon), v  \rangle }_{(c)}+\underbrace{\langle  (v\times B_1^\varepsilon) \cdot\nabla_v(\mu g_1^\varepsilon),v \rangle}_{(d)}=0. 
\]
The first two terms, treated in \cite{BGL1}, give rise to $(a) \rightarrow \partial_t u$ and $(b)\rightarrow u\!\cdot \!\nabla_xu+\nabla_xp-\nu\Delta u$ and here the viscosity $\nu$ depends on the collision operator $L$. See the section 4 in \cite{BGL1} for more detail for the treatment of $(b)$. Note that  $(c)\rightarrow 0$. And for $(d)$, we see that 
\[
(d)\rightarrow \int v  (v\times
B)\cdot\nabla_v({\mu}g) dv=  -J\times B. 
\]
Hence, we obtain the following momentum equation 
\begin{equation}\label{3.7}
\partial_{t}u+u\cdot\nabla u+\nabla p  =\nu\Delta
u+J\times B\,.
\end{equation}
Similarly from the local conservation law of the energy, we can derive the equation for $\theta$. Since 
\[
\begin{split}
&\int \frac{|v|^2}{2}  (v\times
B_1)\cdot\nabla_v({\mu}g_1) dv =- \int (v\times B_1)\cdot v{\mu}g_1dv =0, 
\end{split}
\]
we indeed get the same Fourier equation as in the pure Boltzmann equation 
\begin{equation}\label{3.8}
\partial_{t}\theta+u\cdot\nabla\theta
=\kappa\Delta\theta. 
\end{equation}
Together with the limiting Maxwell's equations, we deduce that  $\rho,u,\theta,J, E,B$  satisfy the incompressible Navier-Stokes-Fourier-Maxwell system: 
\begin{equation}
\begin{split}\label{3.9}
&\bullet \;\rho+\theta=0,\;\; \nabla\cdot u=0\\
&\bullet \;\partial_{t}u+u\cdot\nabla u+\nabla p   =\nu\Delta
u+J\times B\\
&\bullet\; \partial_{t}\theta+u\cdot\nabla\theta
=\kappa\Delta\theta\\
&\bullet\;E+u\times B=\eta 
J\\
&\bullet \;\mu_0\epsilon_0\partial_t E-\nabla\times B=-\mu_0J,\; 
 \partial_tB+\nabla\times
E=0,\;\nabla\cdot B=0
\end{split}
\end{equation}
This system was studied in \cite{Masmoudi10jmpa}. 
 We have established the following result: 

\begin{theorem} Let $[f_1^\varepsilon,g_1^\varepsilon,E_1^\varepsilon,B_1^\varepsilon]$ be the classical solutions to \eqref{rvmbMp}. Assume that  $f_1^\varepsilon$, $g_1^\varepsilon$, $E_1^\varepsilon$, $B_1^\varepsilon$ converge strongly to $f,g,E,B$ as $\varepsilon \rightarrow 0$. Then the limit $f$ has the form $f=\rho+v\cdot u+(\frac{|v|^2}{2}-\frac32)\theta\,.$ Moreover, $\rho,{u},\theta$, $J$, $E,B$ solve the incompressible Navier-Stokes-Fourier-Maxwell system
\eqref{3.9}. 
\end{theorem}

\subsection{Incompressible viscous MHD}

We now discuss the viscous MHD limit. We consider the same scaling as in \eqref{rvmbM} and in addition,  we take $\epsilon_0=\varepsilon$ so that the electric permittivity is sufficiently small: 
\begin{equation}
\begin{split}
&\varepsilon\partial _tF^{\varepsilon}+v\cdot \nabla
_xF^{\varepsilon}+(\varepsilon E^{\varepsilon}+{v\times
B^{\varepsilon}})\cdot \nabla
_v G^{\varepsilon}=\frac{1}{\varepsilon}Q(F^{\varepsilon}, F^{\varepsilon}),\\
&\varepsilon\partial _tG^{\varepsilon}+v\cdot \nabla
_xG^{\varepsilon}+(\frac{E^{\varepsilon}}{\varepsilon}+\frac{v\times
B^{\varepsilon}}{\varepsilon^2})\cdot \nabla
_v F^{\varepsilon}=\frac{1}{\varepsilon}Q(G^{\varepsilon},F^{\varepsilon}),\\
&\mu_0\varepsilon\partial _tE^{\varepsilon}-\nabla \times
B^{\varepsilon}=-{\mu_0} \int_{{\mathbb{ R}}
^3}v \;G^{\varepsilon}\;dv,\;\;\nabla \cdot B^{\varepsilon}=0,\\
&\partial _tB^{\varepsilon}+\nabla \times
E^{\varepsilon}=0,\;\;\nabla \cdot
E^{\varepsilon}=\frac{1}{\varepsilon}\int_{{\mathbb{
R}}^3}G^{\varepsilon}\;dv.\label{rvmb}
\end{split}
\end{equation}
As before, we consider the following perturbation of $F^\varepsilon,G^\varepsilon, E^\varepsilon, B^\varepsilon$ around $[\mu,0,0,0]$: 
\[
F^\varepsilon =\mu (1+\varepsilon f_1^\varepsilon), \;G^\varepsilon=\varepsilon \mu g_1^\varepsilon, \;
E^\varepsilon =\varepsilon E_1^\varepsilon,\; B^\varepsilon =\varepsilon B_1^\varepsilon
\]
and we suppose that there exist the classical solutions $f_1^\varepsilon,\,g_1^\varepsilon,\,E_1^\varepsilon,\, B_1^\varepsilon$ to \eqref{rvmbM} and further assume that as $\varepsilon\rightarrow 0$, $f_1^\varepsilon\rightarrow f$, $g_1^\varepsilon\rightarrow g$, 
$E_1^\varepsilon\rightarrow E$, $B_1^\varepsilon\rightarrow B$. Following the same procedure as in the previous section, we can derive the equations for $\rho,u,\theta,J,E,B$.  The only difference from the previous system \eqref{3.9} lies in the Maxwell's equations. The last equation in \eqref{rvmb} forces $\sigma=0$ when $\varepsilon\rightarrow 0$ and the third equation in \eqref{rvmb} yields the AmpereÕs law $\nabla\times B=\mu_0 J$. We obtain the following theorem:

\begin{theorem} Let $F^\varepsilon =\mu (1+\varepsilon f_1^\varepsilon), \,G^\varepsilon=\varepsilon \mu g_1^\varepsilon, \,
E^\varepsilon =\varepsilon E_1^\varepsilon,\, B^\varepsilon =\varepsilon B_1^\varepsilon$ be the classical solutions to \eqref{rvmb}. Assume that  $f_1^\varepsilon$, $g_1^\varepsilon$, $E_1^\varepsilon$, $B_1^\varepsilon$ converge strongly to $f,g,E,B$ as $\varepsilon \rightarrow 0$. Then the limit $f$ has the form $f=\rho+v\cdot u+(\frac{|v|^2}{2}-\frac32)\theta\,.$ Moreover, $\rho,{u},\theta$, $J$, $E,B$ solve the following incompressible Navier-Stokes-Fourier-MHD equations: 
\begin{equation}
\begin{split}
&\bullet \;\rho+\theta=0,\;\; \nabla\cdot u=0\\
&\bullet\;\partial_{t}u+u\cdot\nabla u+\nabla p   =\nu\Delta
u+J\times B\\
&\bullet\;\partial_{t}\theta+u\cdot\nabla\theta
=\kappa\Delta\theta\\
&\bullet\;E+u\times B=\eta 
J \\
&\bullet\;\nabla\times B=\mu_0J,\;  \partial_tB+\nabla\times
E=0,\;\nabla\cdot B=0 
\end{split}
\end{equation}
\end{theorem}


\subsection{Incompressible inviscid MHD}

To see how the incompressible inviscid MHD system can be derived, we introduce the following scaling to 
\eqref{vmb} by letting $\tilde{t}=\varepsilon^3 t$, $\tilde{x}=\varepsilon^2 x$,  $\tilde{E}= \frac{1}{{\varepsilon}^{3/2}} E$, $\tilde{B}=\frac{1}{\sqrt{\varepsilon}}B$, $\tilde{G}=\frac{1}{{\varepsilon}^{3/2}}G$. By dropping $\tilde{ \text{ }}$, 
\begin{equation}
\begin{split}
&\varepsilon\partial _tF^{\varepsilon}+v\cdot \nabla
_xF^{\varepsilon}+(\varepsilon E^{\varepsilon}+{v\times
B^{\varepsilon}})\cdot \nabla
_v G^{\varepsilon}=\frac{1}{\varepsilon^2}Q(F^{\varepsilon}, F^{\varepsilon}),\\
&\varepsilon\partial _tG^{\varepsilon}+v\cdot \nabla
_xG^{\varepsilon}+\frac{1}{\varepsilon^2}(E^{\varepsilon}+\frac{v\times
B^{\varepsilon}}{\varepsilon})\cdot \nabla
_v F^{\varepsilon}=\frac{1}{\varepsilon^2}Q(G^{\varepsilon},F^{\varepsilon}),\\
&\mu_0\epsilon_0\varepsilon \partial _tE^{\varepsilon}-\frac{1}{\varepsilon^2} \nabla \times
B^{\varepsilon}=-\frac{\mu_0}{\varepsilon^2} \int_{{\mathbb{ R}}
^3}v \;G^{\varepsilon}\;dv,\;\;\nabla \cdot B^{\varepsilon}=0,\\
&\partial _tB^{\varepsilon}+\nabla \times
E^{\varepsilon}=0,\;\;\nabla \cdot
E^{\varepsilon}=\frac{1}{\epsilon_0\varepsilon^2}\int_{{\mathbb{
R}}^3}G^{\varepsilon}\;dv.\label{rvmbIE}
\end{split}
\end{equation}
We consider the following perturbation of $F^\varepsilon,G^\varepsilon, E^\varepsilon, B^\varepsilon$ around $[\mu,0,0,0]$: 
\[
F^\varepsilon =\mu (1+\varepsilon f_1^\varepsilon), \;G^\varepsilon=\varepsilon \mu g_1^\varepsilon, \;
E^\varepsilon =\varepsilon E_1^\varepsilon,\; B^\varepsilon =\varepsilon B_1^\varepsilon. 
\]
We plug this expansion into \eqref{rvmbIE}  to obtain the similar expression to \eqref{rvmbMp}:
\begin{equation}
\begin{split}
&\varepsilon\partial _tf_1^{\varepsilon}+v\cdot \nabla
_xf_1^{\varepsilon}+\varepsilon (\varepsilon E_1^{\varepsilon}+{v\times
B_1^{\varepsilon}})\cdot \frac{1}{\mu}\nabla
_v (\mu g_1^{\varepsilon}) \\
&\quad\quad\quad\quad\quad\quad\quad\quad\quad=\frac{1}{\varepsilon^2 \mu}\{ Q(\mu f_1^\varepsilon, \mu) +Q(\mu,\mu f_1^\varepsilon ) +\varepsilon Q(\mu f_1^\varepsilon,\mu f_1^\varepsilon)\}\\
&\varepsilon\partial _tg_1^{\varepsilon}+v\cdot \nabla
_xg_1^{\varepsilon}+(\frac{E^{\varepsilon}}{\varepsilon^2}+\frac{v\times
B^{\varepsilon}}{\varepsilon^3})\cdot  \frac{1}{\mu}\nabla
_v(\mu+\varepsilon\mu f_1^{\varepsilon}) \\
&\quad\quad\quad\quad\quad\quad\quad\quad\quad=\frac{1}{\varepsilon^2 \mu} \{Q(\mu g_1^\varepsilon , \mu) +\varepsilon Q(\mu g_1^\varepsilon, \mu f_1^{\varepsilon})  \},\\
&\mu_0\epsilon_0\varepsilon \partial _tE_1^{\varepsilon}-\frac{1}{\varepsilon^2}\nabla \times
B_1^{\varepsilon}=-\frac{\mu_0}{\varepsilon^2} \int_{{\mathbb{ R}}
^3}v \mu g_1^{\varepsilon}\;dv,\;\;\nabla \cdot B_1^{\varepsilon}=0,\\
&\partial _tB_1^{\varepsilon}+\nabla \times
E_1^{\varepsilon}=0,\;\;\nabla \cdot
E_1^{\varepsilon}=\frac{1}{\epsilon_0\varepsilon^2}\int_{{\mathbb{
R}}^3}\mu g_1^{\varepsilon}\;dv.\label{rvmbIEp}
\end{split}
\end{equation}
As before, under the convergence assumption of the solutions, we deduce that $f=\rho+v\cdot u+(\frac{|v|^2}{2}-\frac32)\theta\,.$ and also we can  derive the same Ohm's law \eqref{3.6}. We can  follow the same procedure as in the previous sections to derive the equations for $\rho,u,\theta,J,E,B$. The difference is the faster relaxation $\frac{1}{\varepsilon^2}$ than the previous cases and the diffusion for fluid variables doesn't appear at the first order. We can obtain the following theorem. 

\begin{theorem}
Suppose that there exist the classical solutions $f_1^\varepsilon,\,g_1^\varepsilon,\,E_1^\varepsilon,\, B_1^\varepsilon$ to \eqref{rvmbIEp} and further assume that as $\varepsilon\rightarrow 0$, $f_1^\varepsilon\rightarrow f$, $g_1^\varepsilon\rightarrow g$, 
$E_1^\varepsilon\rightarrow E$, $B_1^\varepsilon\rightarrow B$. Then $f=\rho+v\cdot u+(\frac{|v|^2}{2}-\frac32)\theta\,.$ and the limiting equations can be recorded as follows: 
\begin{equation}
\begin{split}
&\bullet \;\rho+\theta=0,\;\; \nabla\cdot u=0\\
&\bullet\;\partial_{t}u+u\cdot\nabla u+\nabla p   =J\times B\\
&\bullet\;\partial_{t}\theta+u\cdot\nabla\theta
=0\\
&\bullet\;E+u\times B=\eta 
J \\
&\bullet\;\nabla\times B=\mu_0J,\;  \partial_tB+\nabla\times
E=0,\;\nabla\cdot B=0 
\end{split}
\end{equation}
which we refer to the incompressible inviscid (Euler) resistive MHD equations. 
\end{theorem}

{We remark that the incompressible, inviscid, ideal MHD equations  are valid when the resistivity $\eta$ is negligible.} 

\subsection{Different masses and Hall effect}

We now look at the case when $m_+\neq m_-$. Moreover, we consider the quasi-neutral regime where $n_+=n_-=1$, we assume the same charges $e_+=e_-=1$, and that we choose $T_0=1$  such that 
\[
\int \mu_+ dv=\int  \mu_- dv=1
\]
where  $\mu_+$ and $\mu_-$ are given in \eqref{GM}. We also assume that  
\[
\int (m_+ \mu_+ + m_-\mu_-)dv=m_++m_-=1. 
\]
We take the same parabolic scaling for $t,x,E,B$ as done in \eqref{rvmb} to \eqref{VMB} except for the scaling of $G$: 
\begin{equation}
\begin{split}
&\varepsilon\partial _tF^\varepsilon_{+}+v\cdot \nabla _xF^\varepsilon_{+}+\frac{1}{m_{+}}(E^\varepsilon+
v\times \frac{B^\varepsilon}{\varepsilon})\cdot \nabla _vF^\varepsilon_{+}=\frac1\varepsilon Q(F^\varepsilon_{+},F^\varepsilon_{+})+\frac1\varepsilon Q(F^\varepsilon_{+},F^\varepsilon_{-}), \\
&\varepsilon\partial _tF^\varepsilon_{-}+v\cdot \nabla _xF^\varepsilon_{-}-\frac{1}{m_{-}}(E^\varepsilon+
v\times\frac{ B^\varepsilon}{\varepsilon})\cdot \nabla _vF^\varepsilon_{-}=\frac1\varepsilon Q(F^\varepsilon_{-},F^\varepsilon_{+})+\frac1\varepsilon Q(F^\varepsilon_{-},F^\varepsilon_{-}). 
\label{rVMBhall} 
\end{split}
\end{equation}
We then consider the following expansion of the solutions near the global Maxwellians $\mu_+$ and $\mu_-$: 
\begin{equation}
\begin{split}\label{exp4}
F^\varepsilon_{+}=\mu_+(1 +\varepsilon f_+^\varepsilon) ,\quad F^\varepsilon_{-}=\mu_- (1+\varepsilon f_-^\varepsilon ), \quad 
E^\varepsilon= \varepsilon E_1^\varepsilon,\quad B^\varepsilon= \varepsilon B_1^\varepsilon\,. 
\end{split}
\end{equation}
The system \eqref{rVMBhall} can be rewritten as 
\begin{equation}
\begin{split}
\varepsilon\partial _tf^\varepsilon_{+}&+v\cdot \nabla _xf^\varepsilon_{+}+\frac{1}{m_{+}\mu_+}(E_1^\varepsilon+
v\times \frac{B_1^\varepsilon}{\varepsilon})\cdot \nabla _v(\mu_++\varepsilon \mu_+ f_+^\varepsilon)\\
&=-\frac{1}{\varepsilon}\mathfrak{L}^1\binom{f^\varepsilon_+}{f^\varepsilon_-}+ \frac{1}{\mu_+}\{ Q(\mu_+f^\varepsilon_{+},\mu_+f^\varepsilon_{+})+Q(\mu_+f^\varepsilon_{+},\mu_-f^\varepsilon_{-})\}, \\
\varepsilon\partial _tf^\varepsilon_{-}&+v\cdot \nabla _xf^\varepsilon_{-}-\frac{1}{m_{-}\mu_-}(E_1^\varepsilon+
v\times\frac{ B_1^\varepsilon}{\varepsilon})\cdot \nabla _v(\mu_-+\varepsilon \mu_-f^\varepsilon_{-})\\
&=-\frac{1}{\varepsilon}\mathfrak{L}^2\binom{f^\varepsilon_+}{f^\varepsilon_-}+\frac{1}{\mu_-}\{ Q(\mu_-f^\varepsilon_{-},\mu_+f^\varepsilon_{+})+ Q(\mu_-f^\varepsilon_{-},\mu_-f^\varepsilon_{-})\},
\label{rVMBhallp} 
\end{split}
\end{equation}
where we have used $\mathfrak{L}^1$ and $\mathfrak{L}^2$ to denote each component of the linearized collision operator $\mathfrak{L}$ given in   \eqref{L}. 
By taking $\varepsilon \rightarrow 0$, we see that 
\[
\mathfrak{L}\binom{f^+}{f^-}=0
\]
namely 
\[
\begin{split}
\binom{f^+(t,x,v)}{f^-(t,x,v)}=n^+(t,x)\binom{1}{0} + n^-(t,x)\binom{0}{1}+ u(t,x)\cdot  \binom{m_+v}{ m_-v} +\theta(t,x) \binom{\frac{m_+|v|^2}{2}-\frac32  }{ \frac{m_-|v|^2}{2}-\frac32}
\end{split}
\]
where $n^+,n^-, u,\theta$ are to be determined.  Then the local conservation laws of mass and momentum yield 
\[
\begin{split}
&\varepsilon \langle \partial _t\binom{f^\varepsilon_{+}}{f^\varepsilon_-}, \binom{1}{0}\rangle_M + \nabla_x\cdot \langle\binom{f^\varepsilon_{+}}{f^\varepsilon_-}, \binom{v}{0}\rangle_M=0, \\
&\varepsilon \langle \partial _t\binom{f^\varepsilon_{+}}{f^\varepsilon_-}, \binom{0}{1}\rangle_M + \nabla_x\cdot \langle\binom{f^\varepsilon_{+}}{f^\varepsilon_-}, \binom{0}{v}\rangle_M=0, \\
&\varepsilon \langle \partial _t \binom{f^\varepsilon_{+}}{f^\varepsilon_-},\binom{m_+ v}{m_-v} \rangle_M  + \nabla_x\!\cdot\!\langle\binom{f^\varepsilon_{+}}{f^\varepsilon_-}, \binom{m_+v\otimes v}{m_-v\otimes v}\rangle_M \\
&\quad\quad\quad+ \langle (\varepsilon  E_1^\varepsilon+v\times {B_1^\varepsilon})\!\cdot\!\nabla_v(\mu_+f_+^\varepsilon-\mu_-f_-^\varepsilon) ,v\rangle =0. 
\end{split}
\]
As $\varepsilon\rightarrow 0$, since $ \langle (v\times B_1)\!\cdot\!\nabla_v(\mu_+f_+-\mu_-f_-) ,v\rangle = -
 \langle v,\mu_+f_+-\mu_-f_-\rangle \times B_1 =0\times B_1=0$, we deduce that 
\[
\begin{split}
 \nabla_x\!\cdot\! \langle\binom{f_{+}}{f_-}, \binom{v}{0}\rangle_M=0, \; \nabla_x\!\cdot\! \langle\binom{f_{+}}{f_-}, \binom{v}{0}\rangle_M=0  \quad\Longrightarrow\quad   \nabla_x\!\cdot\! u=0 \\
\nabla_x\!\cdot\!\langle\binom{f_{+}}{f_-}, \binom{m_+v\otimes v}{m_-v\otimes v}\rangle_M=0 \quad\Longrightarrow\quad \nabla_x(\rho+\theta)=0
\end{split}
\]
where $\rho=n_++n_-$. We remark that 
\begin{equation}\label{J_0}
J_0\equiv\int v (\mu_+f_+ -\mu_-f_-)dv=0.
\end{equation}
By assuming that 
\begin{equation}\label{JJ}
\lim_{\varepsilon\rightarrow 0}\frac{1}{\varepsilon}J^\varepsilon=J_1
\end{equation}
one can derive the MHD type equations containing $J_1\times B_1$ as before.

To discuss the Ohm's law, we write out the dynamics of the charge density $G^\varepsilon:=\mu_+f^\varepsilon_+-\mu_-f^\varepsilon_-$ from \eqref{rVMBhallp}: 
\begin{equation}
\begin{split}\label{ohmg}
&\varepsilon \partial_t G^\varepsilon +v\!\cdot\!\nabla_xG^\varepsilon + (E_1^\varepsilon+
v\times \frac{B_1^\varepsilon}{\varepsilon})\cdot \nabla _v(\frac{\mu_+}{m_+}+\frac{\mu_-}{m_-} +\varepsilon  (\frac{\mu_+ f_+^\varepsilon}{m_+}+\frac{\mu_- f_-^\varepsilon}{m_-}))\\
&=-\frac{1}{\varepsilon}\left(\mu_+\mathfrak{L}^1\binom{f^\varepsilon_+}{f^\varepsilon_-}  -\mu_-\mathfrak{L}^2\binom{f^\varepsilon_+}{f^\varepsilon_-}\right)\\
&\quad+\{Q(\mu_+f^\varepsilon_{+},\mu_+f^\varepsilon_{+})+Q(\mu_+f^\varepsilon_{+},\mu_-f^\varepsilon_{-}) - Q(\mu_-f^\varepsilon_{-},\mu_+f^\varepsilon_{+})- Q(\mu_-f^\varepsilon_{-},\mu_-f^\varepsilon_{-}) \}
\end{split}
\end{equation}
Now projecting \eqref{ohmg} onto $v$, we obtain the dynamics of $J^\varepsilon=\langle v G^\varepsilon \rangle$: 
\begin{equation}\label{ohmlaw}
\begin{split}
&\varepsilon \partial_t J^\varepsilon +\nabla_x\!\cdot\! \langle v\otimes v G^\varepsilon \rangle - (\frac{1}{m_+}+\frac{1}{m_-} +\varepsilon (\frac{n_+^\varepsilon}{m_+}+\frac{n_-^\varepsilon}{m_-} ) ) E^\varepsilon_1\\&- \underbrace{(\frac{\langle v,\mu_+f_+^\varepsilon\rangle}{m_+}+\frac{\langle v,\mu_-f_-^\varepsilon\rangle}{m_-}) \times B^\varepsilon_1 }_{(\ast)}
=- \frac{1}{\varepsilon}\eta J^\varepsilon  + \langle v, Q(\mu_+f^\varepsilon_{+},\mu_-f^\varepsilon_{-}) - Q(\mu_-f^\varepsilon_{-},\mu_+f^\varepsilon_{+}) \rangle 
\end{split}
\end{equation}
To derive the Ohm's law, we split the last term $(\ast)$ in the left-hand-side of \eqref{ohmlaw}  into two parts:   
\begin{equation}
\begin{split}
(\ast)&=  \frac{\langle v,m_+\mu_+f_+^\varepsilon+m_-\mu_-f^\varepsilon_-\rangle}{m_+m_-} \times B^\varepsilon_1 +  (\frac{1}{m_+}-\frac{1}{m_-}) \langle v,\mu_+f_+^\varepsilon-\mu_-f^\varepsilon_-\rangle \times B^\varepsilon_1\\
 &=\frac{1}{m_+m_-} u^\varepsilon \times B^\varepsilon_1 +  \underbrace{(\frac{1}{m_+}-\frac{1}{m_-}) J^\varepsilon\times B_1^\varepsilon}_{(\ast\ast)}. 
\end{split}\label{ast}
\end{equation}
By \eqref{J_0}, $(\ast\ast)\rightarrow 0$. Hence the Ohm's law for $J_1\equiv \lim_{\varepsilon\rightarrow 0}\frac{1}{\varepsilon}J^\varepsilon$  is given by 
\[
\eta J_1= \frac{1}{m_+m_-} (E_1 +u\times B_1 ) +\nabla_x P + \mathcal{C}
\]
in the limit $\varepsilon\rightarrow 0$. Here $[E_1,B_1]=\lim_{\varepsilon\rightarrow 0}[ E_1^\varepsilon,B_1^\varepsilon]$,  $u=\lim_{\varepsilon\rightarrow 0}\langle v,m_+\mu_+f_+^\varepsilon+m_-\mu_-f^\varepsilon_-\rangle $, $\nabla_xP=-\nabla_x\!\cdot\!\langle v\otimes v G_0\rangle$ where $G_0=\lim_{\varepsilon\rightarrow 0} G^\varepsilon$, and  the pressure $P$ may depend on $\frac{n^-}{m_-}-\frac{n^+}{m_+}$  and $(\frac{1}{m_-}-\frac{1}{m_+})\theta$.  Notice that we didn't have the pressure term in the previous simplified cases because we started with the strong quasi-neutral scaling where $G_0=0$.  The last term $\mathcal{C}$ coming from the last term in \eqref{ohmlaw} is given by $\langle v, Q(\mu_+ f_+, \mu_- f_-)-Q(\mu_-f_-,\mu_+f_+) \rangle$. In fact, $f_+$ and $f_-$ are purely hydrodynamic described by $n^\pm$, $u$, $\theta$ and hence $\mathcal{C}$ may depend on $n^\pm$, $u$, $\theta$ quadratically.  Notice that $\mathcal{C}$ vanishes when $m_+=m_-$. 

Now we briefly discuss the Hall effect in the above setting. As in the compressible case, the Hall effect is not seen in the first order limit. The Hall effect is represented by the Hall term $J_1\times B_1$, which should be seen from ${(\ast\ast)}$ in \eqref{ast}.  Since $\lim_{\varepsilon\rightarrow 0}J^\varepsilon =0$ by \eqref{J_0}, $(\ast\ast)$ is 0 at the leading order. If we approximate $J^\varepsilon$ by $\varepsilon J_1$ by using \eqref{JJ}, we see that  $(\ast\ast)=\varepsilon (\frac{1}{m_+}-\frac{1}{m_-})J_1\times B_1+o(\varepsilon)$.  Therefore, we deduce  that the Hall term $J_1\times B_1$  in \eqref{ohmlaw} can be captured at the next order by including  $O(\varepsilon)$ terms in the spirit of Chapman-Enskog type expansions. It is interesting to point out that $(\ast\ast)$ is always 0 for $m_+=m_-$ and thus it will not affect the Ohm's law at all; for instance, see Section 3.1 and 3.2.  This confirms that the Hall effect is a {\it{two-fluid effect}}.  In the case of $m_+\neq m_-$,  the general Ohm's law \eqref{ohmlaw} is much more complicated in that not only the Hall effect but also other higher order effects might appear. This motivates us to study the higher order expansion which will be investigated in a separate article.

\section{Conclusion}

We would like to conclude the paper with few remarks: 
 
1) In this paper, we tried to give a formal derivation in the spirit of \cite{BGL1} to 
MHD type models starting form the Vlasov-Maxwell-Boltzmann system. 
For compressible models, we were able to derive the 
$1\frac12$ Euler-Maxwell system. It seems to us that the classical 
two fluid  Euler-Maxwell system can not be derived from the model we started 
from and we plan to come back to this problem in a forthcoming work

2) We were also able to derive Ohm's law and to compute the resistivity 
from the Boltzmann kernel. This can be compared to the formulas giving 
the viscosity 
$\nu$ and heat diffusivity  $\kappa$ in \cite{BGL1}. 
It is worth pointing out that in our derivation the Hall effect which yields 
a correction to the Ohm's law can only be seen  as a 
higher order correction.  This is similar to the compressible Navier-Stokes 
system which is a correction to the  compressible  Euler system.

3) There are few other important scaling parameters that we did not take 
into account in these derivation. One of them is the mass ratio between 
the electron and the ions. Usually it reflects the fact that the 
electrons attain their equilibrium much faster than the ions.

4)  As observed by Grad \cite{Grad}, Ohm's law is also valid through an entirely different mechanism, the so-called gyro-oscillation. We will leave this (magnetic) gyro-effect with or without collisions for future study.

5) After most of this work was done, we learnt that D.~Arsenio and L. Saint-Raymond  \cite{AS11}
are also studying similar limits starting form the renormalized solutions 
of DiPerna and Lions.

\

\textbf{Acknowledgements.} The authors would like to thank the referees for helpful comments and suggestions, which have improved the presentation of the paper.  Part of work was done during the authors' visit to ICERM, Brown University in the fall of 2011 and they thank the institute for the support.

\end{document}